\documentclass[11pt,reqno]{amsart}
\usepackage{enumerate, latexsym, amsmath, amsfonts, amssymb, amsthm, color}
\def\pmod #1{\ ({\rm{mod}}\ #1)}
\def\Z{\Bbb Z}

\def\Q{\Bbb Q}

\def\l{\left}
\def\r{\right}
\def\bg{\bigg}
\def\({\bg(}
\def\){\bg)}
\def\t{\text}
\def\f{\frac}

\def\ls{\leqslant}

\def\ve{\varepsilon}

\def\eq{\equiv}

\def\da{\delta}

\def\Proof{\noindent{\it Proof}}
\def\Ack{\medskip\noindent {\bf Acknowledgment}}
\theoremstyle{plain}
\newtheorem{theorem}{Theorem}

\newtheorem{lemma}{Lemma}
\newtheorem{corollary}{Corollary}

\theoremstyle{definition}

\theoremstyle{remark}

\allowdisplaybreaks

 \vspace{4mm}

\begin{document}

\hbox{Czechslovak Math. J. 73 (2023), no.\,3, 971--978.}
\medskip

\title
[{The tangent function and power residues modulo primes}]
{The tangent function and \\ power residues modulo primes}

\author
[Zhi-Wei Sun] {Zhi-Wei Sun}

\address{Department of Mathematics, Nanjing
University, Nanjing 210093, People's Republic of China}
\email{zwsun@nju.edu.cn}

\keywords{Power residues modulo primes, the tangent function, identity.
\newline \indent The research has been supported
by the National Natural Science Foundation of China (grant 11971222).}

\begin{abstract} Let $p$ be an odd prime, and let $a$ be an integer not divisible by $p$.
When $m$ is a positive integer with $p\eq1\pmod{2m}$ and $2$ is an $m$th power residue modulo $p$,
we determine the value of the product $\prod_{k\in R_m(p)}(1+\tan\pi\frac{ak}p)$, where
$$R_m(p)=\{0<k<p:\ k\in\mathbb Z\ \text{is an}\ m\text{th power residue modulo}\ p\}.$$
In particular, if $p=x^2+64y^2$ with $x,y\in\mathbb Z$, then
$$\prod_{k\in R_4(p)}\left(1+\tan\pi\frac {ak}p\right)=(-1)^{y}(-2)^{(p-1)/8}.$$
\end{abstract}
\maketitle

\section{Introduction}
\setcounter{lemma}{0}
\setcounter{theorem}{0}
\setcounter{corollary}{0}
\setcounter{remark}{0}
\setcounter{equation}{0}

It is well known that the function $\tan \pi x$ has period $1$. For any positive odd number $n$
and complex number $x$ with $x-1/2\not\in\Z$, Sun \cite[Lemma 2.1]{PMD} proved that
$$\prod_{r=0}^{n-1}\l(1+\tan\pi\f{x+r}n\r)=\l(\f 2n\r)2^{(n-1)/2}\l(1+\l(\f{-1}n\r)\tan\pi x\r),$$
where $(\f{\cdot}n)$ is the Jacobi symbol. In particular, for any odd prime $p$ and integer $a\not\eq0\pmod p$ we have
$$\prod_{k=1}^{p-1}\l(1+\tan\pi\f{ak}p\r)=\prod_{r=0}^{p-1}\l(1+\tan\pi\f rp\r)=\l(\f 2p\r)2^{(p-1)/2}.$$

Let $p$ be an odd prime. Then
$$1^2,\ 2^2,\ \ldots,\ \l(\f{p-1}2\r)^2$$
modulo $p$ give all the $(p-1)/2$ quadratic residues modulo $p$.
Sun \cite[Theorem 1.4]{PMD} determined the value of the product
$\prod_{k=1}^{(p-1)/2}(1+\tan\pi\f{ak^2}p)$
for any integer $a$ not divisible by $p$; in particular,
\begin{align*}&\prod_{k=1}^{(p-1)/2}\l(1+\tan\pi\f{ak^2}p\r)
\\=\ &\begin{cases}(-1)^{|\{1\ls k<\f p4:\ (\f kp)=1\}|}2^{(p-1)/4}&\t{if}\ p\eq1\pmod8,
\\(-1)^{|\{1\ls k<\f p4:\ (\f kp)=-1\}|}2^{(p-1)/4}(\f ap)\ve_p^{-3(\f ap)h(p)}&\t{if}\ p\eq5\pmod8,
\end{cases}
\end{align*}
where $(\f{\cdot}p)$ is the Legendre symbol, and $\ve_p$ and $h(p)$ are the fundamental unit and the class number of the real quadratic field
$\Q(\sqrt p)$ respectively.

Let $m\in\Z^+=\{1,2,3,\ldots\}$, and let $p$ be a prime with $p\eq1\pmod{m}$.
If $a\in\Z$ is not divisible by $p$, and $x^m\eq a\pmod p$ for an integer $x$,
then $a$ is called an $m$th power residue modulo $p$.
The set
\begin{equation}\label{R_m(p)}
R_m(p)=\{k\in\{1,\ldots,p-1\}:\ k\ \t{is an}\ m\t{th power residue modulo}\ p\}
\end{equation}
has cardinality $(p-1)/m$, and $\{k+p\Z:\ k\in R_m(p)\}$ is a subgroup of the multiplicative group
$\{k+p\Z:\ k=1,\ldots,p-1\}$.
For an integer $a\not\eq0\pmod p$, the $m$th power residue symbol
$(\f ap)_m$ is a unique $m$th root $\zeta$ of unity such that
$$a^{(p-1)/m}\eq\zeta\pmod p$$
in the ring of all algebraic integers. (Note that a primitive root $g$ modulo $p$
has order $p-1$ which is a multiple of $m$.) In particular,
$$\l(\f{-1}p\r)_m=(-1)^{(p-1)/m}.$$

Let $p$ be a prime with $p\eq1\pmod{2m}$, where $m\in\Z^+$.
Note that $p-1\in R_m(p)$ since $(-1)^{(p-1)/m}=1$.
If $2\in R_m(p)$, then $-2=(-1)\times2$ is an $m$th power residue modulo $p$, hence
$$\l(\f{-2}p\r)_{2m}=\begin{cases}1&\t{if}\ -2\ \t{is a}\ 2m\t{-th power residue modulo}\ p,
\\-1&\t{otherwise},\end{cases}$$
and
\begin{equation}\label{-2/p}\l(\f{-2}p\r)_{2m}^m=\l(\f{-2}p\r)
\end{equation}
since
$$\l(\f{-2}p\r)_{2m}^m\eq\l((-2)^{(p-1)/(2m)}\r)^m=(-2)^{(p-1)/2}\eq\l(\f{-2}p\r)\pmod p.$$

Now we state our main theorem.

\begin{theorem} \label{Main} Let $m\in\Z^+$, and let $p$ be a prime with $p\eq1\pmod{2m}$.
Suppose that $2$ is an $m$th power residue modulo $p$. For any integer $a$ not divisible by $p$, we have
\begin{equation}\label{tan-xy} \prod_{k\in R_m(p)}\l(1+\tan\pi\f {ak}p\r)=\l(\f{-2}p\r)_{2m}(-2)^{(p-1)/(2m)}=\l(\f 2p\r)_{2m}2^{(p-1)/(2m)}.
\end{equation}
\end{theorem}

To prove Theorem \ref{Main}, we need the following auxiliary result.

\begin{theorem}\label{Main'} Let $m$ be a positive integer, and let $p$ be a prime
with $p\eq1\pmod{2m}$. Suppose that $2$ is an $m$th power residue modulo $p$. For any integer $a\not\eq0\pmod p$, we have
\begin{equation}\label{G(i)}\prod_{k\in R_m(p)}(i-e^{2\pi iak/p})=\l(\f{-2}p\r)_{2m}i^{(p-1)/(2m)}
\end{equation}
and
\begin{equation}\label{G(i)+}\prod_{k\in R_m(p)}(i+e^{2\pi iak/p})=\l(\f{2}p\r)_{2m}i^{(p-1)/(2m)}.
\end{equation}
\end{theorem}

Let $p$ be an odd prime with $p\eq 1\pmod m$, where $m$ is $3$ or $4$. Then there are unique $x,y\in\Z^+$ such that $p=x^2+my^2$ (cf. \cite[pp.\,7-12]{Cox}). It is well known that $2\in R_m(p)$ if and only if $p=x^2+m(my)^2$ for some $x,y\in\Z^+$ (cf. Prop. 9.6.2 of \cite[p.\,119]{IR}
and Exer. 26 of \cite[p.\,64]{IR}).

Theorem \ref{Main} with $m=3$ has the following consequence.

\begin{corollary} \label{Cor1.1} Let $p=x^2+27y^2$ be a prime with $x,y\in\Z^+$. For any integer $a\not\eq0\pmod p$, we have
\begin{equation}\label{tan-xy} \prod_{k\in R_3(p)}\l(1+\tan\pi\f {ak}p\r)=(-1)^{xy/2}(-2)^{(p-1)/6}.
\end{equation}
\end{corollary}

From Theorem \ref{Main} in the case $m=4$, we can deduce the following result.

\begin{corollary} \label{Cor1.2} Let $p=x^2+64y^2$ be a prime with $x,y\in\Z^+$. For any integer $a\not\eq0\pmod p$, we have
\begin{equation}\label{tan-4} \prod_{k\in R_4(p)}\l(1+\tan\pi\f {ak}p\r)=(-1)^{y}(-2)^{(p-1)/8}.
\end{equation}
\end{corollary}

We will prove Theorems 1.1-1.2 in the next section,
 and deduce Corollaries 1.1-1.2 in Section 3.

\section{Proofs of Theorems \ref{Main} and \ref{Main'}}
\setcounter{lemma}{0}
\setcounter{theorem}{0}
\setcounter{corollary}{0}
\setcounter{remark}{0}
\setcounter{equation}{0}

\medskip

\begin{lemma}\label{Lem2.1} Let $m$ be a positive integer, and let $p$ be a prime with $p\eq1\pmod{2m}$.
Then we have
$$\sum_{k\in R_m(p)}k =\f{p(p-1)}{2m}.$$
\end{lemma}
\Proof. Note that $-1$ is an $m$th power residue modulo $p$ since $(p-1)/m$ is even.
For $k\in\{1,\ldots,p-1\}$, clearly $p-k\in R_m(p)$ if and only if $k\in R_m(p)$.
Thus
$$2\sum_{k\in R_m(p)}k=\sum_{k\in R_m(p)}(k+(p-k))=p\times|R_m(p)|=\f{p(p-1)}{m}.$$
This ends the proof of Lemma 2.1. \qed

\medskip

\noindent{\it Proof of Theorem \ref{Main'}}.
Let
$$c:=\prod_{k\in R_m(p)}\l(i-e^{2\pi iak/p}\r).$$
As $k\in\Z$ is an $m$th power residue modulo $p$ if and only if $-k$
is an $m$th power residue modulo $p$, we also have
$$c=\prod_{k\in R_m(p)}\l(i-e^{2\pi ia(-k)/p}\r).$$
Thus
\begin{align*}c^2&=\prod_{k\in R_m(p)}\l(i-e^{2\pi iak/p}\r)\l(i-e^{-2\pi iak/p}\r)
\\&=\prod_{k\in R_m(p)}\l(i^2+1-i\l(e^{2\pi iak/p}+e^{-2\pi iak/p}\r)\r)
\\&=(-i)^{|R_m(p)|}\prod_{k\in R_m(p)}\l(e^{2\pi iak/p}+e^{-2\pi iak/p}\r)
\\&=(-i)^{(p-1)/m}\prod_{k\in R_m(p)}e^{-2\pi iak/p}\l(1+e^{4\pi i ak/p}\r)
\\&=(-1)^{(p-1)/(2m)}e^{-2\pi i\sum_{k\in R_m(p)}ak/p}\prod_{k\in R_m(p)}\f{1-e^{2\pi ia(4k)/p}}{1-e^{2\pi ia(2k)/p}}.
\end{align*}
Note that
$$e^{-2\pi i\sum_{k\in R_m(p)}ak/p}=e^{-2\pi ia(p-1)/(2m)}=1$$
by Lemma \ref{Lem2.1}. As $2$ is an $m$th power residue modulo $p$, we also have
$$\prod_{k\in R_m(p)}\l(1-e^{2\pi i ak/p}\r)=\prod_{k\in R_m(p)}\l(1-e^{2\pi i a(2k)/p}\r)
=\prod_{k\in R_m(p)}\l(1-e^{2\pi i a(4k)/p}\r).$$
Combining the above, we see that
$$c^2=(-1)^{(p-1)/(2m)}\times1\times1=(-1)^{(p-1)/(2m)}.$$

Write $c=\da i^{(p-1)/(2m)}$ with $\da\in\{\pm1\}$. In the ring of all algebraic integers, we have
\begin{align*}c^p&=\prod_{k\in R_m(p)}(i-e^{2\pi iak/p})^p
\\&\eq\prod_{k\in R_m(p)}(i^p-1)=(i^p-1)^{(p-1)/m}
\\&=((i^p-1)^2)^{(p-1)/(2m)}=(-2i^p)^{(p-1)/(2m)}\pmod p.
\end{align*}
Thus
$$\da i^{p(p-1)/(2m)}=c^p\eq(-2)^{(p-1)/(2m)}i^{p(p-1)/(2m)}\pmod p$$
and hence
$$\da \eq(-2)^{(p-1)/(2m)}\eq\l(\f{-2}p\r)_{2m}\pmod p.$$
Therefore $\da=(\f{-2}p)_{2m}$ and hence \eqref{G(i)} holds.

Taking conjugates of both sides of \eqref{G(i)}, we get
$$\prod_{k\in R_m(p)}(-i-e^{-2\pi i ak/p})=\l(\f{-2}p\r)_{2m}(-i)^{(p-1)/(2m)}$$
and hence
$$(-1)^{(p-1)/m}\prod_{k\in R_m(p)}(i+e^{2\pi ia(p-k)/p})=\l(\f{-2}p\r)_{2m}\l(\f{-1}p\r)_{2m}i^{(p-1)/(2m)}.$$
This is equivalent to \eqref{G(i)+} since $\{p-k:\ k\in R_m(p)\}=R_m(p)$.

In view of the above, we have completed the proof of Theorem 1.2. \qed

\medskip

\noindent{\it Proof of Theorem \ref{Main}}.
For any $k\in\Z$, we have
\begin{align*}1+\tan\pi\f kp&=1+\f{\sin \pi k/p}{\cos\pi {k}/p}
=1+\f{(e^{ i\pi k/p}-e^{-i\pi k/p})/(2i)}{(e^{ i\pi k/p}+e^{-i\pi k/p})/2}
\\&=1-i\f{e^{2\pi i k/p}+1-2}{e^{2\pi ik/p}+1}=1-i+\f{2i}{e^{2\pi i k/p}+1}
\\&=(1-i)\l(1+\f{i-1}{e^{2\pi ik/p}+1}\r)=(1-i)\f{e^{2\pi ik/p}+i}{e^{2\pi i k/p}-i^2}
\\&=\f{i-1}{i-e^{2\pi ik/p}}\times\f{e^{2\pi i(2k)/p}-i^2}{e^{2\pi ik/p}-i^2}.
\end{align*}
Therefore
\begin{equation}\label{mid} \prod_{k\in R_m(p)}\l(1+\tan\pi\f {ak}p\r)=\f{(i-1)^{|R_m(p)|}}{\prod_{k\in R_m(p)}(i-e^{2\pi iak/p})}.
\end{equation}

Recall \eqref{G(i)} and note that
$$(i-1)^{|R_m(p)|}=((i-1)^2)^{(p-1)/(2m)}=(-2i)^{(p-1)/(2m)}.$$
So \eqref{mid} yields that
\begin{align*}\prod_{k\in R_m(p)}\l(1+\tan\pi\f {ak}p\r)&=\f{(-2i)^{(p-1)/(2m)}}{(\f{-2}p)_{2m}i^{(p-1)/(2m)}}
\\&=\l(\f{-2}p\r)_{2m}(-2)^{(p-1)/(2m)}
\\&=\l(\f{2}p\r)_{2m}2^{(p-1)/(2m)}.
\end{align*}
This concludes our proof of Theorem 1.1. \qed

\section{Proofs of Corollaries 1.1-1.2}
\setcounter{lemma}{0}
\setcounter{theorem}{0}
\setcounter{corollary}{0}
\setcounter{remark}{0}
\setcounter{equation}{0}

\begin{lemma}\label{Lem3.1} For any prime $p=x^2+27y^2$ with $x,y\in\Z^+$, we have
\begin{equation}\label{xy} \l(\f{-2}p\r)=(-1)^{xy/2},
\end{equation}
\end{lemma}
\Proof. Clearly $p\eq1\pmod6$ and $x\not\eq y\pmod 2$ since $p=x^2+27y^2$.
Note that \eqref{xy} has the equivalent form:
\begin{equation}\label{4|xy} 4\mid xy\iff p\eq1,3\pmod 8.
\end{equation}
\medskip

{\it Case} 1. $x$ is odd and $y$ is even.
\medskip

In this case,
$$p=x^2+27y^2\eq1+3y^2=1+12\l(\f y2\r)^2\eq1+4\l(\f y2\r)^2\pmod 8$$
and hence
$$p\eq1,3\pmod 8\iff p\eq1\pmod8\iff 2\mid \f y2\iff 4\mid y\iff 4\mid xy.$$

\medskip

{\it Case} 2. $x$ is even and $y$ is odd.
\medskip

In this case,
$$p=x^2+27y^2\eq x^2+3y^2=4\l(\f x2\r)^2+3\pmod 8$$
and hence
$$p\eq1,3\pmod 8\iff p\eq3\pmod8\iff 2\mid \f x2\iff 4\mid x\iff 4\mid xy.$$

In view of the above, we have completed the proof of Lemma \ref{Lem3.1}. \qed

\medskip

\noindent{\it Proof of Corollary \ref{Cor1.1}}. As $p=x^2+27y^2$, we see that $p\eq1\pmod6$
and $2$ is a cubic residue modulo $p$. By Lemma 3.1 and \eqref{-2/p} with $m=3$, we have
$$\l(\f{-2}p\r)_{6}=\l(\f{-2}p\r)=(-1)^{xy/2}.$$
Combining this with Theorem \ref{Main} in the case $m=3$, we immediately obtain the desired
\eqref{tan-xy}. \qed

\medskip

\noindent{\it Proof of Corollary \ref{Cor1.2}}. As $p=x^2+64y^2$, we see that $p\eq1\pmod8$
and $2$ is a quartic residue modulo $p$. By Theorem 7.5.7 or Corollary 7.5.8 of \cite[pp.\,227-228]{BEW}, we have
$$\l(\f{-2}p\r)_{8}=(-1)^{y}.$$
Combining this with Theorem \ref{Main} in the case $m=4$, we immediately obtain the desired
\eqref{tan-4}. \qed

\medskip
\Ack. The author would like to thank the referee for helpful comments.


\medskip


\begin{thebibliography}{99}

\bibitem{BEW} B. C. Berndt, R. J. Evans and K. S. Williams, Gauss and Jacobi Sums,
John Wiley \& Sons, 1998.

\bibitem{Cox} D. A. Cox,  Primes of the Form $x^2+ny^2$: Fermat, Class Field Theory and Complex Multiplication, John Wiley \& Sons, Inc., New York, 1989.


\bibitem{IR} K. Ireland and M. Rosen, {\it A Classical Introduction to Modern Number Theory},
2nd Edition, Grad. Texts. in Math. 84, Springer, New York, 1990.


\bibitem{PMD} Z.-W. Sun, {\it Trigonometric identities and quadratic residues},
Publ. Math. Debrecen {\bf 102} (2023), 111--138.


\end{thebibliography}
 \end{document}